# Dirichlet eta and beta functions at negative integer arguments: Exact results from anti-limits

Kamal Bhattacharyya*

Harish-Chandra Research Institute, Chhatnag Road, Jhusi, Allahabad 211019, India

**Abstract**

A route to evaluate exact sums represented by Dirichlet eta and beta functions, both of which are alternating and divergent at negative integer arguments, is advocated. It rests on precise polynomial extrapolations and stands as a generalization of an early endeavor on lattice sums. Apart from conferring a physical meaning to anti-limits, the scheme advanced here is direct, independent and computationally appealing. A new interpretation of summability is also gained.

**Keywords:** Dirichlet eta and beta functions; divergent saw-tooth sequence; anti-limit; lattice sums

**PACS:** 02.30.Gp; 02.30.Lt

**MSC:** 11M41; 30B50; 40A05

*E-mail: kamalbhattacharyya@hri.res.in; pchemkb@gmail.com

# 1. Introduction

Saw-tooth sequences emerge from alternating series. The limit point $S_\infty$ of a convergent sequence $\{S_j\}$ of partial sums thus defines the actual sum of the parent series. Such sequences diverge when the parent series diverges. However, assigning a 'sum' to a specific divergent series is not always a straightforward task. Hardy [1] discussed at length on the *sense* in which such 'sums' are to be interpreted. Dirichlet eta and beta functions at negative integer arguments provide considerable insights [1] in the context of divergent *alternating* series. On the other hand, Riemann zeta function belongs to the other category, diverging *monotonically* at negative integer arguments. Some of the members of these two types are mutually linked too [1]. Here, however, we shall be concerned with the former functions only for which a simple geometrical interpretation exists, as detailed below.

The problems associated with ionic lattice sums [2-8] provide nice practical testing grounds to understand and treat divergences of saw-tooth sequences, deciphering specifically the meaning of *sums* and their functional relevance. In the course of our early explorations on ionic lattice sums [6-8], even with isovalent impurities and defects [8], we adopted Padé-type sequence acceleration schemes on various occasions to gain considerable success. A few other approaches were also in vogue (see, *e.g.*, references quoted in [6-8]). Nonetheless, meanwhile [7] we also suggested an extrapolation scheme (ES) and observed that the *anti-limit* of a divergent saw-tooth sequence could be *physical* as well, so much so that even an approximate sequential scheme [7] would be quite rewarding, within limits. Indeed, during our early exploration, it was also evident that *exact* linear extrapolation (LE) [7] agrees with some of Hardy's axioms [1] that stand as foundations of handling any divergent series. Here, we provide a generalization of the approach.

A preliminary survey of the relevant literature reveals that, while a lot of works focused attention on eta and beta functions (see, *e.g.*, [9-11] and references quoted therein), they are chiefly restricted to positive arguments. In contrast, zeta functions received attention in respect of studies in both computations and related properties for positive arguments [12-15], taming along with the divergence at negative integer arguments [16, 17]. Therefore, a simple, direct and exact methodology for evaluations of the two functions at issue in the concerned domain seems demanding.

In view of the above remarks, purpose of the present communication is fourfold.

First, we wish to establish that the ES *is* the right choice to tackle both Dirichlet eta and beta functions at negative integer arguments. Second, these series are generally amenable to *exact* polynomial extrapolation (PE), and hence the results are also *exact*. Further, our recipe obeys Hardy's axioms [18] too. Third, this endeavor brings to light the physical role of the *anti-limit*. Fourth, the divergent eta and beta functions may now be justifiably called as 'PE summable', in addition to the existing ones [1].

## 2. A preamble to the anti-limit

Consider the partial sums $S_1$, $S_2$, $S_3$, ... of an alternating series that form a sawtooth sequence $\{S_M\}$. The sequence tends to the correct answer $S_\infty$ when the parent series is convergent. Otherwise, the limit point $S_\infty$ does not exist. For a divergent alternating series, we have found [7] it expedient to split the corresponding sequence $\{S_M\}$ into two parts as *odd* $\{S_o\}$ [*e.g.*, $S_1$, $S_3$, $S_5$, ...] and *even* $\{S_e\}$ [*e.g.*, $S_2$, $S_4$, $S_6$, ...]. In simple situations, members of $\{S_o\}$ and $\{S_e\}$ fall *exactly* on two distinctly different straight lines defined by $P_o(x)$ and $P_e(x)$ where $x$ is the *sequence number*. They meet at some point(s) $X$ towards the left, such that $P_o(X) = P_e(X) = P(X) = S_X$, the value at that point. Thus, in such divergent situations, the anti-limits $S_X$ replace the conventional $S_\infty$. This forms the basis of LE in the context of ES.

A few examples will nicely clarify the situation. Let us first choose $\eta(-1)$, where

$$\eta(-1) = 1 - 2 + 3 - 4 + \dots. \tag{1}$$

Figure 1 shows how the result ¼ emerges. Here, $\{S_o\}$ [*e.g.*, 1, 2, 3, ...] yield $P_o(x) = x/2 + ½$ (blue line) while $\{S_e\}$ [*e.g.*, -1, -2, -3, ...] satisfy $P_e(x) = -x/2$ (red line). Solving, we find $X$ and $P(X)$. The LE is *exact*, so also the outcome, as noted in the figure. A similar situation prevails for the $\beta(-1)$ series, expressed as

$$\beta(-1) = 1 - 3 + 5 - 7 + \dots. \tag{2}$$

Now, we have $P_o(x) = x$ and $P_e(x) = -x$, leading to $X = 0$ and $S_X = \beta(-1) = 0$. As the next illustration of the LE, we explore how an answer to $\eta(0)$ may be found, where

$$\eta(0) = 1 - 1 + 1 - 1 + \dots. \tag{3}$$

To this end, we add up (1) and (3), obtaining

$$\eta(-1) + \eta(0) = 2 - 3 + 4 - 5 + \dots \tag{4}$$

Employing now the LE to (4), one gets $P_o(x) = x/2 + 3/2$ and $P_e(x) = -x/2$; hence $X = -3/2$ and $S_X = ¾$ follow. Thus, with known $\eta(-1)$, it turns out that $\eta(0) = ½$. In effect, this result is true for $\beta(0)$ too, since it parallels expansion (3). Alternatively, one may consider the sequence for

$\beta(0) + \beta(-1)$. Another quick execution concerns the value of $\zeta(0)$ that goes as

$$\zeta(0) = 1+1+1+1+... \tag{5}$$

Adding up (1) and (5), we notice

$$\eta(-1) + \zeta(0) = 2-1+4-3+.... \tag{6}$$

From (6), one obtains $P_o(x) = 3x/2 + \frac{1}{2}$ and $P_e(x) = x/2$, implying $X = -\frac{1}{2}$ and $S_X = -\frac{1}{4}$. With the already known value of $\eta(-1) = \frac{1}{4}$, (6) yields $\zeta(0) = -\frac{1}{2}$. The consistency and benefit of the LE are thus obvious.

The above idea had earlier been exploited to handle a few *practical* lattice-sum sequences that are divergent in character. There we adopted an approximate *sequential* LE scheme [7] with reasonable success, but considerable physical flavor. Surprisingly, we have now unearthed a more general *exact* ES via the PE.

### 3. The procedure

In a nutshell, our procedure will run as follows: (i) Check if the series is alternating and divergent, and then proceed. (ii) Compute members of the set $\{S_o\}$ up to a certain number, depending on the choice of *s*. (iii) Fit the numbers via a polynomial $P_o(x)$ such that it exactly reproduces the input at odd integer *x*. (iv) Ensure that the polynomial does not change even if the input data or the degree of the polynomial is increased. (v) Repeat steps (ii) to (iv) for computed members of the set $\{S_e\}$ to find $P_e(x)$. This polynomial should exactly reproduce the input at even integer *x*. (vi) Solve subsequently for $x = X$ at which $P_o(x) = P_e(x)$ is satisfied. (vii) Evaluate $P_o(X) = P_e(X)$ to get the answer.

Two more points need to be emphasized. First, the recipe applies only to a *class* of alternating divergent series for which *exact* polynomial representations for $P_o(x)$ and $P_e(x)$ exist (*e.g.*, power series stand as notable exceptions), yielding, in general, $X \leq 2$ (anti-limit). This is understandable because, in a convergent case, we would have $X = \infty$ (limit point). So, our strategy is not *regular*. Second, the point(s) $X$ need not *always* be real; indeed, this condition may sometimes be relaxed (see Sec. 5.3).

### 4. Results and discussion

It will now be seen how the PE becomes the *natural candidate* to tame both $\eta(s)$ and $\beta(s)$ series for any negative integer *s*. The policy, results (down to $s = -10$) and associated discussion are expounded below.

### 4.1. The Dirichlet eta function

The Dirichlet eta function is defined as

$$\eta(s) = \sum_{k=1}^{\infty} \frac{(-1)^{k-1}}{k^s}. \tag{7}$$

The series is alternating and converges *absolutely* for any positive integer $s \geq 2$. Problems start from $s = 1$. Specifically, one finds $\eta(1) = \ln 2$ by *inspection*. It is *conditionally* convergent. We already noted how $\eta(0)$ and $\eta(-1)$ are obtained. However, the situation becomes more frustrating as we gradually go downwards. The LE is helpless too.

As prescribed earlier, we now seek whether the even and odd sequences separately possess *exact* polynomial representations $P_e^\eta(x)$ and $P_o^\eta(x)$ for a given *s*-value. If so, an extrapolation to the left ($X \leq 2$) will reveal that there is *at least* one *point of intersection*. It might then be identified that $P^\eta(X) \equiv \eta(s)$ for that specific *s*. The physical implication is clear; the parent alternating sequence $\{S_M\}$ appears to diverge from the common point(s) $X$.

In Figure 2, we show the case of $\eta(-3)$. Here, the fittings of $P_o^\eta(x), P_e^\eta(x)$ by cubic polynomials are *exact*, but one needs larger polynomials for higher $|s|$. Table 1 displays such results. Figure 2 also shows that the lines meet near zero and the value is slightly negative. The actual result is given in Table 1. Interestingly, while the simplest LE ($s = -1$) possesses a *single X* that has already appeared in Figure 1, a careful scrutiny uncovers *multiple X*-values for *any s* < -1. But, a great relief is that they *always* yield the *same $\eta(s)$*. Figure 3 shows this assuring behavior for the cases under consideration here, starting from $\eta(-2)$ to $\eta(-10)$, exposing *all possible X*-values. However, our interpretation rests on the *first intersection point*, as the lines approach from the right side.

A few interesting properties of the aforesaid polynomials [*cf.* Table 1] may now be in order: (i) The degree of a polynomial is dictated by $-s$. (ii) $P_e^\eta(x)$ is devoid of any constant term. (iii) If *s* is *even*, there is no constant term in $P_o^\eta(x)$ as well. Thus, $P_e^\eta(x) = -P_o^\eta(x)$, and hence the result becomes zero [*cf.* cases of $\eta(-2), \eta(-4), ...,$ etc.]. Here, $X = 0$ is always *a common intersection point* that may be easily checked. (iv) $P_e^\eta(x)$ is always expressible as $P_e^\eta(x) = -[P_o^\eta(x) - k_o^\eta]$ wherever a positive constant $k_o^\eta$ appears in $P_o^\eta(x)$. This constant term additionally signifies that $\eta(s) = k_o^\eta/2$. (v) For *odd s*, a *common intersection point* is $X = -\frac{1}{2}$, as Figure 3 shows. (vi) Structure of these polynomials is also notable. The highest power of $x$ (equal to $|s|$) is followed by the immediately lower one, and then the still lower ones appear alternately. This also explains why no constant term turns up for *even s*.

### *4.2. The Dirichlet beta function*

The Dirichlet beta function (also known as $L(s)$, especially in Hardy [1]) is defined as

$$\beta(s) = \sum_{k=1}^{\infty} \frac{(-1)^{k-1}}{(2k-1)^s}. \tag{8}$$

Once again, LE yield here the exact answers for $\beta(-1)$ and $\beta(0) + \beta(-1)$. These have been noted already. At any other $s < -1$, fittings for $P_o^\beta(x), P_e^\beta(x)$ are *exact*, though, as before, larger polynomials are needed for higher $|s|$. Behaviors of such polynomials *around zero* are presented in Figure 4. Multiple $X$-values are seen here as well, but all offering the *same* $P^\beta(X) = \beta(s)$. This is again unambiguous, and hence comforting. Table 2 displays the various pairs of polynomials over the range $-1 \geq s \geq -10$.

As before, we may summarize now a few properties of $P_o^\beta(x)$ and $P_e^\beta(x)$ point-wise: (i) The degree of a polynomial is dictated by $-s$. (ii) $P_e^\beta(x)$ is devoid of any constant term. (iii) If $s$ is *odd*, there is no constant term in $P_o^\beta(x)$ as well. Thus, $P_e^\beta(x) = -P_o^\beta(x)$, and hence the result becomes *zero* [*cf.* cases of $\beta(-3), \beta(-5), \beta(-7),$ etc.]. Thus, $X = 0$ is a common point of intersection in these situations. (iv) It is always possible to express the even polynomial $P_e^\beta(x)$ as $P_e^\beta(x) = -[P_o^\beta(x) - k_o^\beta]$ wherever a positive constant $k_o^\beta$ appears in $P_o^\beta(x)$. This constant term additionally signifies that $\beta(s) = k_o^\beta/2$. (v) When $s$ is *even*, $X = \frac{1}{2}$ is a common point, as Figure 4 reveals. (vi) Structure of these polynomials is also notable. The highest power of $x$ is followed by the lower ones alternately. As a result, no constant term survives for *odd s*.

## 5. Further remarks

### *5.1. General*

In Sec. 3, we remarked briefly about the approach to follow. To be wise after the events, one may be eager to simplify it further. For $\eta(s)$, one should continue with steps (ii) to (iv) first. The degree of the fitting polynomial is equal to $|s|$. Half the constant term in $P_o(x)$ will be the answer, and it is *non-zero* only for *odd s*. All these characteristics may be verified.

The structures of the polynomials also act as yardsticks. Alongside, one may further check that the corresponding $P_e(x)$ will generate members of the $\{S_e\}$. The polynomials are *odd* about $x = -½$ when $s$ is odd, and *even* when $s$ is even. In case of $\beta(s)$, the evaluation of $P_o(x)$ and the value would follow the same route. Half the constant term in $P_o(x)$ will again be the answer, but here it is *non-zero* only for *even s*. Once these traits match, $P_e(x)$ is easily constructed to inspect the emergence of correct $\{S_e\}$. Additionally, here the polynomials are *odd* about $x = 0$ when $s$ is odd, and *even* when $s$ is even. In order to convince ourselves further, we tabulate the results for two more *still lower s*-values below. Table 3 endorses the above assertions beyond any doubt. Moreover, the following properties of the characteristic polynomials may serve as final checks for a general reader:

$s = \text{even} : P_o^\eta(x=0) = P_e^\eta(x=0) = 0;\ P_o^\eta(x=-1) = P_e^\eta(x=-1) = 0.$

$s = \text{odd} : P_e^\eta(x=0) = 0;\ P_o^\eta(x=-1) = 0.$

$s = \text{even} : P_e^\beta(x=0) = 0.$

$s = \text{odd} : P_e^\beta(x=0) = P_o^\beta(x=0) = 0.$

### 5.2. Special choices

Let us now turn attention to some additional returns. The usual practice of estimating $\eta(s)$ rests either on the relation

$$\eta(s) = \left(1 - 2^{1-s}\right)\zeta(s) \qquad (9)$$

that requires known $\zeta(s)$, or a type of symmetry relation of the form

$$\eta(-s) = 2s \frac{1 - 2^{-(1+s)}}{(1 - 2^{-s})\pi^{s+1}} \sin\left(\frac{\pi s}{2}\right) \Gamma(s)\, \eta(s+1). \qquad (10)$$

In either case, we have an advantage. Knowing from Table 3 the result

$$\eta(-19) = 1 - 2^{19} + 3^{19} - 4^{19} + \ldots = -\tfrac{221930581}{8}, \qquad (11)$$

one obtains from (9) the following outcome:

$$\zeta(-19) = 1 + 2^{19} + 3^{19} + 4^{19} + \ldots = \tfrac{174611}{6600}. \qquad (12)$$

In other words, values of certain *monotonically* divergent series are obtainable indirectly from our results. Similarly, using relations like (11), it is possible to retrieve $\eta$-values for positive integer arguments from (10). The latter series are alternating but *convergent*. For the beta functions, likewise, we have a relation that reads as

$$\beta(1-s) = \left(\tfrac{\pi}{2}\right)^{-s} \sin\left(\frac{\pi s}{2}\right) \Gamma(s)\beta(s). \tag{13}$$

Obtaining the value of $\beta(-20)$ from Table 3, *e.g.*,

$$\beta(-20) = 1 - 3^{20} + 5^{20} - 7^{20} + \ldots = \tfrac{370371188237525}{2}, \tag{14}$$

one gets $\beta(21)$ from (13), which is the sum of again an alternating but *convergent* series. Note that the converse process is usually adopted in (10) or (13). One computes the right side to obtain the left. But, this process has a catch. Unless one finds the *exact* value of the convergent series appearing at the right, which is computationally unrealistic, it is impossible to infer that the answer for the left side will be a *rational number*! In our direct scheme, no such problems arise. In this way, the present endeavor accomplishes a greater good.

### *5.3. Hardy's axioms*

Hardy's axioms [18] are as follows:

$$(A): \sum_n a_n = s \Rightarrow \sum_n \mu a_n = \mu s.$$
$$(B): \sum_n a_n = s, \sum_n b_n = t \Rightarrow \sum_n (a_n + b_n) = s + t.$$
$$(C): a_1 + a_2 + a_3 + \ldots = s \Rightarrow v + a_1 + a_2 + a_3 + \ldots = v + s.$$

We now outline how they apply to the present context. (i) Axiom (*A*) is followed here because both the odd and even polynomials, along with the answer, are multiplied by $\mu$; the *X*-value(s) remains unaltered, nevertheless. Axiom (*B*) is obeyed likewise, though the symmetry [*cf.* Sec. 5.1] and *X*-value(s) of the resultant polynomial change. For instance, the result of a linear combination like $\beta(-2) + \eta(-3)$ yields the new sequence 2, -15, 37, …, from which we recover the corresponding polynomials, *e.g.*,

$$\begin{aligned}S_o &: 2, 37, 130, \ldots \Rightarrow P_o(x) = -\tfrac{5}{4} + \tfrac{1}{2}x^3 + \tfrac{11}{4}x^2; \\ S_e &: -15, -76, -207, \ldots \Rightarrow P_e(x) = -P_o(x) - \tfrac{5}{4}.\end{aligned} \tag{15}$$

The *X*-values change, but the correct result [-5/8] is regained. As regards axiom (*C*), we checked that the constant $v$ is added to both the polynomials. Consequently, the qualitative character of the intersecting polynomials remains unaltered, apart from a shift along the ordinate by the amount $v$. The result accordingly obeys the third axiom.

A somewhat more general question is, given two polynomials as $P_o(x) = x^2$ and $P_e(x)$ = -1, does the anti-limit strategy apply? It turns out to be relevant because there is no point of intersection along the *real* line. However, they do become equal at *imaginary* $X = \pm i$, with the desired answer -1. We explain the outcome now in terms of Hardy's axioms. The series in this case reads as $1 - 2 + 10 - 10 + 26 - 26 + \ldots$ . Suppose, $\upsilon = 10 + 26 + \ldots$; then, by Hardy's (A), $-\upsilon = -10 - 26 - \ldots$ . We apply now axiom (B) to obtain $0 = 10 - 10 + 26 - 26 + \ldots$ . Choosing $v = 1 - 2$ in (C), we clearly see the correspondence with our answer. The choice $P_o(x) = x^2 + x$ and $P_e(x) = -1$ likewise yield two complex $X$-values for the series $2 - 3 + 13 - 13 + 31 - 31 + \ldots$ . But, here too the answer -1 again agrees with the standard result along similar lines as above.

## 6. Conclusion

We have pursued a direct ES to decipher via the anti-limits how the values of two classes of alternating divergent series could be zero, positive and negative. After the initial success of the LE [7], here we have adopted a more general PE to achieve the end. We happily note that divergent eta and beta functions admit *exact* PEs, and hence anti-limits yield *exact answers* for integers. Thus, special status of *integers* is obvious here. No *transformation* has been applied anywhere. Indeed, *exact* results from *computations* are *rarely* found for any divergent series. Hence, such startling and direct observations allow us to put forward a new definition, *viz.*, the 'PE summability'. In other words, the eta and beta functions at negative integer arguments are PE *summable*. Cases wherever *exact* PEs *follow* for divergent saw-tooth sequences, our definition should apply. For example, no negative *non-integer* argument in the present context admits *exact* PE. The same is true of alternating but divergent power series. We have also noted that anti-limits may even exist in the *complex* domain [see Sec. 5.3]. Thus, along with its simplicity, generality of the prescription is noteworthy too. Above all, a different interpretation of the sum of a class of divergent, alternating series is gained.

**Acknowledgement**

It is a pleasure to acknowledge the Harish-Chandra Research Institute, Prayagraj, for providing all the necessary facilities to pursue the work.

Table 1. Characteristic polynomials for eta functions at different *s*-values down to *s* = -10 with the results found from the anti-limit *X*.

| s | Polynomial | $\eta(s)$ |
|---|---|---|
| -1 | $P_o^\eta(x) = \frac{1}{2}x + \frac{1}{2}; P_e^\eta(x) = -[P_o^\eta(x) - \frac{1}{2}]$ | $\frac{1}{4}$ |
| -2 | $P_o^\eta(x) = \frac{1}{2}x^2 + \frac{1}{2}x; P_e^\eta(x) = -[P_o^\eta(x)]$ | 0 |
| -3 | $P_o^\eta(x) = \frac{1}{2}x^3 + \frac{3}{4}x^2 - \frac{1}{4}; P_e^\eta(x) = -[P_o^\eta(x) + \frac{1}{4}]$ | $-\frac{1}{8}$ |
| -4 | $P_o^\eta(x) = \frac{1}{2}x^4 + x^3 - \frac{1}{2}x; P_e^\eta(x) = -[P_o^\eta(x)]$ | 0 |
| -5 | $P_o^\eta(x) = \frac{1}{2}x^5 + \frac{5}{4}x^4 - \frac{5}{4}x^2 + \frac{1}{2}; P_e^\eta(x) = -[P_o^\eta(x) - \frac{1}{2}]$ | $\frac{1}{4}$ |
| -6 | $P_o^\eta(x) = \frac{1}{2}x^6 + \frac{3}{2}x^5 - \frac{5}{2}x^3 + \frac{3}{2}x; P_e^\eta(x) = -[P_o^\eta(x)]$ | 0 |
| -7 | $P_o^\eta(x) = \frac{1}{2}x^7 + \frac{7}{4}x^6 - \frac{35}{8}x^4 + \frac{21}{4}x^2 - \frac{17}{8}; P_e^\eta(x) = -[P_o^\eta(x) + \frac{17}{8}]$ | $-\frac{17}{16}$ |
| -8 | $P_o^\eta(x) = \frac{1}{2}x^8 + 2x^7 - 7x^5 + 14x^3 - \frac{17}{2}x; P_e^\eta(x) = -[P_o^\eta(x)]$ | 0 |
| -9 | $P_o^\eta(x) = \frac{1}{2}x^9 + \frac{9}{4}x^8 - \frac{21}{2}x^6 + \frac{63}{2}x^4 - \frac{153}{4}x^2 + \frac{31}{2}; P_e^\eta(x) = -[P_o^\eta(x) - \frac{31}{2}]$ | $\frac{31}{4}$ |
| -10 | $P_o^\eta(x) = \frac{1}{2}x^{10} + \frac{5}{2}x^9 - 15x^7 + 63x^5 - \frac{255}{2}x^3 + \frac{155}{2}x; P_e^\eta(x) = -[P_o^\eta(x)]$ | 0 |

Table 2. Characteristic polynomials for beta functions at different *s*-values with the results found from the anti-limit *X*.

| s | Polynomial | $\beta(s)$ |
|---|---|---|
| -1 | $P_o^\beta(x) = x; P_e^\beta(x) = -P_o^\beta(x)$ | 0 |
| -2 | $P_o^\beta(x) = 2x^2 - 1; P_e^\beta(x) = -[P_o^\beta(x) + 1]$ | $-\frac{1}{2}$ |
| -3 | $P_o^\beta(x) = 4x^3 - 3x; P_e^\beta(x) = -P_o^\beta(x)$ | 0 |
| -4 | $P_o^\beta(x) = 8x^4 - 12x^2 + 5; P_e^\beta(x) = -[P_o^\beta(x) - 5]$ | $\frac{5}{2}$ |
| -5 | $P_o^\beta(x) = 16x^5 - 40x^3 + 25x; P_e^\beta(x) = -P_o^\beta(x)$ | 0 |
| -6 | $P_o^\beta(x) = 32x^6 - 120x^4 + 150x^2 - 61; P_e^\beta(x) = -[P_o^\beta(x) + 61]$ | $-\frac{61}{2}$ |
| -7 | $P_o^\beta(x) = 64x^7 - 336x^5 + 7000x^3 - 427x; P_e^\beta(x) = -P_o^\beta(x)$ | 0 |
| -8 | $P_o^\beta(x) = 128x^8 - 896x^6 + 2800x^4 - 3416x^2 + 1385; P_e^\beta(x) = -[P_o^\beta(x) - 1385]$ | $\frac{1385}{2}$ |
| -9 | $P_o^\beta(x) = 256x^9 - 2304x^7 + 10080x^5 - 20496x^3 + 12465x; P_e^\beta(x) = -P_o^\beta(x)$ | 0 |
| -10 | $P_o^\beta(x) = 512x^{10} - 5760x^8 + 33600x^6 - 102480x^4 + 124650x^2 - 50521;$ $P_e^\beta(x) = -[P_o^\beta(x) + 50521]$ | $-\frac{50521}{2}$ |

Table 3. Characteristic polynomials for eta and beta functions at two low *s*-values.

| s | Polynomial | Value |
|---|---|---|
| -19 | $P_o^\eta(x) = \frac{1}{2}x^{19} + \frac{19}{4}x^{18} - \frac{969}{8}x^{16} + 2907x^{14} - \frac{214149}{4}x^{12} + \frac{1431859}{2}x^{10}$ $- \frac{26113581}{4}x^8 + 37041963x^6 - \frac{900752361}{8}x^4 + \frac{547591761}{4}x^2 - \frac{221930581}{4}$ | $-\frac{221930581}{8}$ |
| | $P_o^\beta(x) = 262144x^{19} - 11206656x^{17} + 317521920x^{15} - 6779092992x^{13}$ $+107193415680x^{11} - 1194759408128x^9 + 8715963060480x^7$ $-37090711793088x^5 + 75161501074020x^3 - 45692713833379x$ | 0 |
| -20 | $P_o^\eta(x) = \frac{1}{2}x^{20} + 5x^{19} - \frac{285}{2}x^{17} + 3876x^{15} - 82365x^{13} + 1301690x^{11}$ $-14507545x^9 + 105834180x^7 - \frac{900752361}{2}x^5 + 912652935x^3 - \frac{1109652905}{2}x$ | 0 |
| | $P_o^\beta(x) = 524288x^{20} - 24903680x^{18} + 793804800x^{16} - 19368837120x^{14}$ $+357311385600x^{12} - 4779037632512x^{10} + 43579815302400x^8 -$ $247271411953920x^6 + 751615010740200x^4 - 913854276667580x^2$ $+370371188237525$ | $\frac{370371188237525}{2}$ |

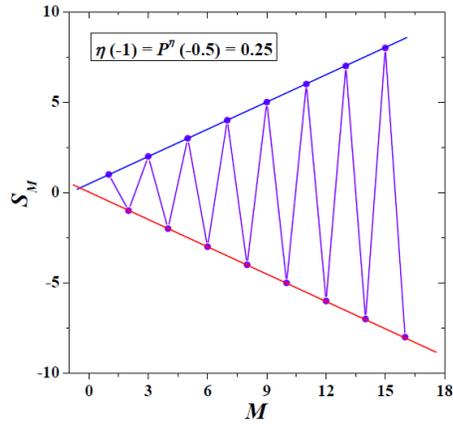

**Figure 1. LE for the case of *η*(-1), yielding the exact result (blue line for *odd* members and red one for *even* members of the sequence).**

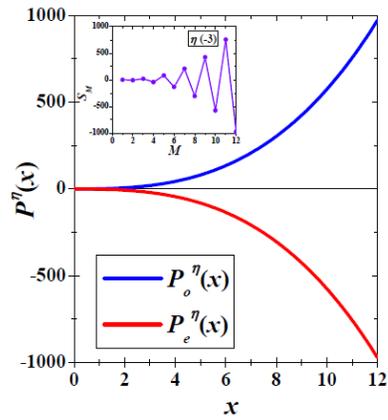

**Figure 2. PE for the sample case of *η*(-3). The inset features the partial sum sequence. The blue line stands for $P_o(x)$ and the red one for $P_e(x)$.**

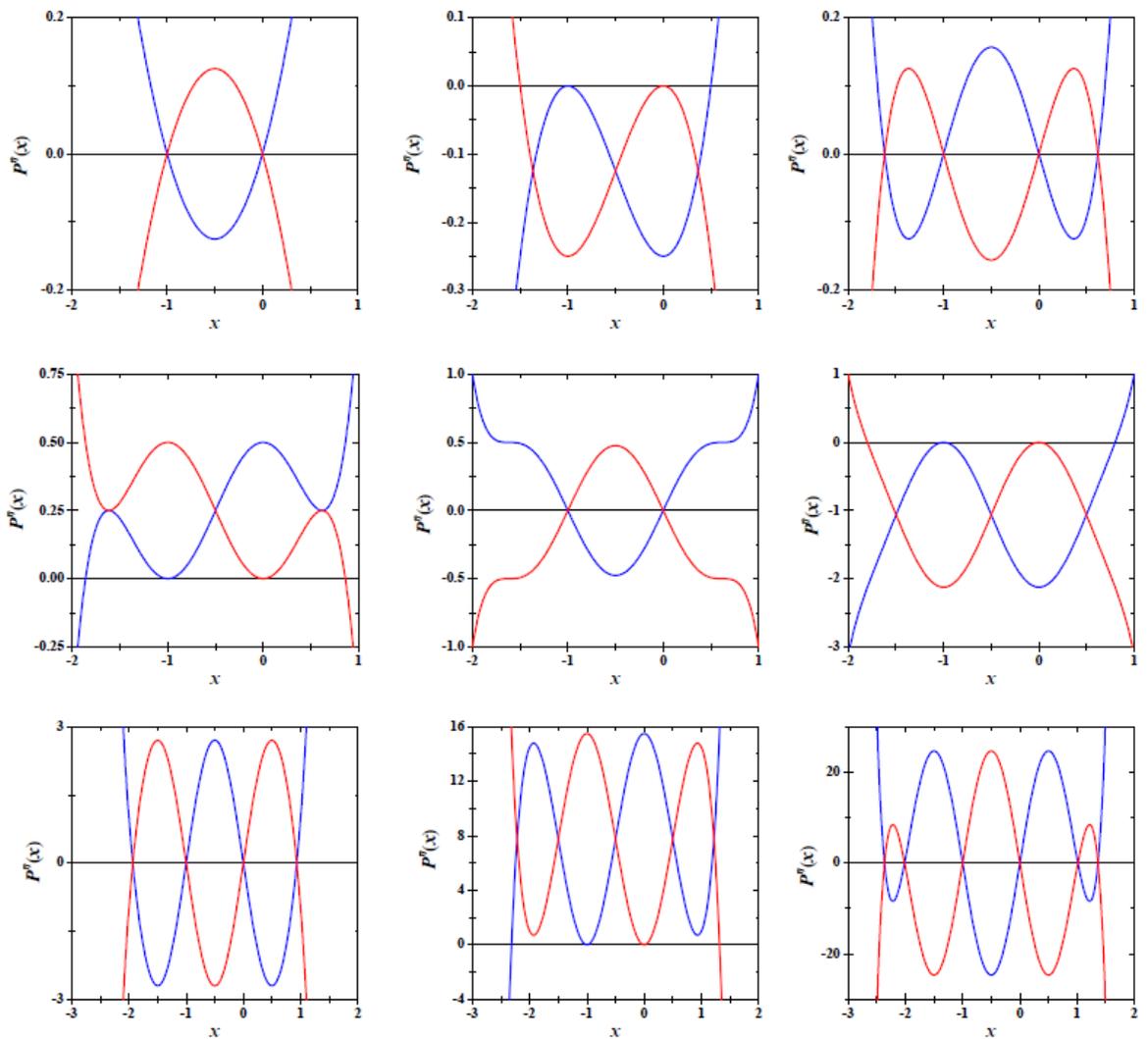

**Figure 3. Nature of $P^{\eta}(x)$ around $X$ for $s$-values from -2 to -10 . Blue and red lines refer respectively to odd and even polynomials. The first row shows cases of $s$ = -2, -3 and -4 as one proceeds towards the right; the last row similarly refers to cases of $s$ = -8, -9 and -10.**

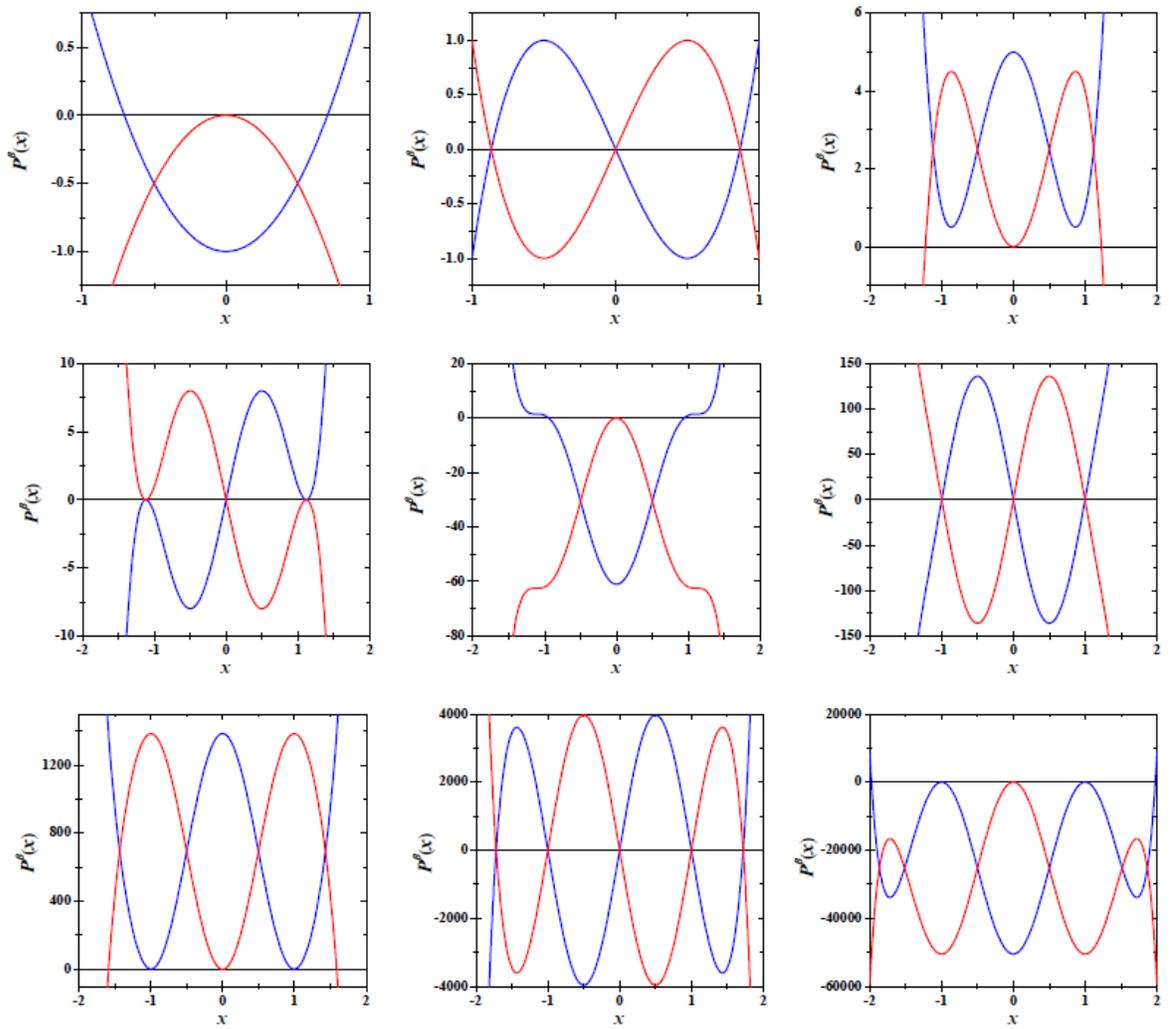

**Figure 4.** Nature of $P^\beta(x)$ around $X$ for $s$ from -2 to -10. All other features are same as that of Figure 3.